\documentclass[a4paper, 10pt, oneside, onecolumn,sort&compress]{elsarticle}
\usepackage{amssymb}
\usepackage{mathrsfs}
\usepackage{amsfonts}
\usepackage[title]{appendix}
\usepackage{syntonly}
\newcommand{\lap}{\mbox{$\bigtriangleup$}}
\newcommand{\ra}{{\mbox{$\rightarrow$}}}
\newtheorem{thm}{Theorem}[section]

\newdefinition{remark}{Remark}[section]
\newdefinition{example}{Example}[section]
\newdefinition{proposition}{Proposition}[section]
\newdefinition{definition}{Definition}[section]
\newproof{proof}{Proof}
\newproof{pot}{Proof of Theorem \ref{the1.1}}
\usepackage{amssymb,amsmath}

\numberwithin{equation}{section}
\journal{Journal of Mathematical Analysis and Applications}

\begin{document}
\begin{frontmatter}

\title{\textbf{Symmetry and Nonexistence of Solutions to Semilinear Equations Involving the Fractional Laplacian on $\mathbb{R}^n$ and $\mathbb{R}^n_+$}\tnoteref{mytitlenote}}
\tnotetext[mytitlenote]{This work is supported by the National Natural Science Foundation of China (No. 11271299) and the Natural Science Basic Research Plan in Shaanxi Province of China (No. 2016JM1023).}
\author{Lizhi Zhang\corref{cor1}}
\ead{azhanglz@163.com}
\author{Yongzhong Wang\corref{cor2}}
\ead{wangyz@nwpu.edu.cn}

\cortext[cor2]{Corresponding author. Tel.: +86 029 88431660.}

\address{Department of Applied Mathematics, Northwestern Polytechnical
University, Xi'an, Shaanxi, 710129, P. R. China}

\begin{abstract}
In this paper, we investigate the following semilinear equations involving the fractional Laplacian
\begin{equation*}
(-\lap)^{\alpha/2} u(x)=f(u),
\end{equation*}
on $\mathbb{R}^n$ and $\mathbb{R}^n_+$ respectively, $0<\alpha<2$. Applying a direct method of moving planes for the fractional Laplacian, we prove symmetry and nonexistence of positive solutions on $\mathbb{R}^n$ and $\mathbb{R}^n_+$ under mild conditions on $f$.
\end{abstract}


\end{frontmatter}


\section{Introduction}
\label{Section 1}

Symmetry and nonexistence properties are very important and useful in studying semilinear elliptic equations. For example, these properties played an essential role in deriving a priori bounds for solutions in \cite{APP}, \cite{FLN}, \cite{GS1}, \cite{HL}, and they were used to obtain uniqueness of solutions in \cite{CC}, \cite{KK}, \cite{MK}, \cite{MS}, \cite{PS}. There are many other applications.

Gidas and Spruck \cite{GS} proved that there is no nontrivial $C^2$ solution to the problem
\begin{equation}
\left\{\begin{array}{ll}
-\lap u(x)=u^p, & \qquad \text{in}~\mathbb{R}^n, \\
u(x)\geq0, & \qquad \text{in}~\mathbb{R}^n,
\end{array}\right.
\label{i01}\end{equation}
if $n\geq3$ and $1<p<\frac{n+2}{n-2}$. Gidas, Ni and Nirenberg in \cite{GNN} pointed out that in the critical case $p=\frac{n+2}{n-2}$ the only solutions of (\ref{i01}) with some decay at infinity are of the form
$$
u(x)=\frac{[n(n-2)\lambda^2]^{\frac{n-2}{4}}}{(\lambda^2+|x-x^0|^2)^{\frac{n-2}{2}}}
$$
for some $\lambda>0$ and $x^0\in \mathbb{R}^n$. This elegant result was also established by Caffarelli, Gidas and Spruck \cite{CGS} without any condition at infinity. Later Chen and Li \cite{CL} found a simpler proof of the result based on the Kelvin transform and the Alexandrov-Serrin moving plane method.

On the upper half Euclidean space
$$
\mathbb{R}^n_+=\{x=(x_1,\cdots,x_n):x_n>0\},
$$
Gidas and Spruck \cite{GS} showed that there are no nontrivial solutions to the Dirichlet problem
\begin{equation}
\left\{\begin{array}{lll}
-\lap u(x)=u^p, & \qquad \text{in}~\mathbb{R}^n_+, \\
u(x)\geq0, & \qquad \text{in}~\mathbb{R}^n_+, \\
u(x)=0, & \qquad \text{on}~\partial\mathbb{R}^n_+,
\end{array}\right.
\label{i02}\end{equation}
where $n\geq3$ and $1<p\leq\frac{n+2}{n-2}$.

These results were generalized to the problems with more general nonlinearity $g(u)$. Bianchi \cite{BG} concluded that if $u\in C^2(\mathbb{R}^n)$ is a positive solution to
\begin{equation}
\lap u+g(u)=0, ~~~\text{in} ~\mathbb{R}^n,
\label{izlz}
\end{equation}
where $g(t)\geq0$ is locally Lipschitz in $t>0$ and $g(t)/t^{\frac{n+2}{n-2}}$
is non-increasing, then either $u$ is a constant or
\begin{equation}
u(x)=\frac{k}{(|x-x^0|^2+h^2)^{(n-2)/2}}
\label{i03}
\end{equation}
for some positive constants $k,~h$ and some point $x^0\in\mathbb{R}^n$. Damascelli and Gladiali \cite{DG} investigated the weak solutions $u\in W^{1,2}_{loc}(\mathbb{R}^n)\cap C^0(\mathbb{R}^n)$ of problem (\ref{izlz}) under the slightly different conditions and obtained the same conclusion.

Using the same technique, Damascelli and Gladiali \cite{DG} also considered the weak solutions to the corresponding Dirichlet problem in $\mathbb{R}^n_+$
\begin{equation}
\left\{\begin{array}{lll}
-\lap u(x)=g(u), & \qquad \text{in}~\mathbb{R}^n_+, \\
u(x)\geq0, & \qquad \text{in}~\mathbb{R}^n_+, \\
u(x)=0, & \qquad \text{on}~\partial\mathbb{R}^n_+,
\end{array}\right.
\label{i04}\end{equation}
where $g:[0,+\infty)\rightarrow\mathbb{R}$ is a continuous function satisfying

$(i)$ $g(t)/t^{\frac{n+2}{n-2}}$ \emph{is nonincreasing in} $(0,+\infty)$;

$(ii)$ $g^+(t)/t$ \emph{is bounded for} $t\rightarrow0$;

$(iii)$ $g(s)>0$ \emph{for every} $s>0$, $\liminf_{s\rightarrow\infty}g(s)>0.$\\
They derived that the solution $u$ depends only on $x_n$ under $(i)$ and $(ii)$, and $u$ is trivial under $(i),~(ii)$ and $(iii)$.

In this paper, we are concern with symmetry and nonexistence of positive solutions for semilinear elliptic equations involving the \textit{fractional Laplacian} on $\mathbb{R}^n$ and $\mathbb{R}^n_+$. Let us begin with the definition of the fractional Laplacian.

The fractional Laplacian in $\mathbb{R}^n$ is a nonlocal pseudo-differential operator with the form
\begin{equation}
(-\lap)^{\alpha/2}u(x)=C_{n,\alpha}PV\int_{\mathbb{R}^n}\frac{u(x)-u(z)}{|x-z|^{n+\alpha}}dz,
\label{pp}
\end{equation}
where $0<\alpha<2$ and $PV$ stands for the Cauchy principle value.

In recent years, the fractional Laplacian has been frequently used to model diverse physical phenomena, such as the turbulence, water waves, the anomalous diffusion, quasi-geostrophic flows, relativistic quantum mechanics of stars and molecular dynamics (see \cite{BoG} \cite{Co} \cite{CaV} \cite{TZ} and the references therein). It also has various applications in probability and finance (see \cite{A} \cite{Be} \cite{CT}). In particular, the fractional Laplacian can be understood as the infinitestmal generator of a stable L\'{e}vy process \cite{Be}. We refer the readers to Di Nezza, Palatucci and Valdinoci's survey paper \cite{NPV}.

The operator in (\ref{pp}) is well defined in the Schwartz space ${\cal S}$ of rapidly decreasing $C^\infty$ functions on $\mathbb{R}^n$. In the space ${\cal S}$, it can be equivalently defined
by the Fourier transform
$$
\widehat{(-\lap)^{\alpha/2}}u(\xi)=|\xi|^\alpha\hat{u}(\xi),
$$
where $\hat{u}$ is the Fourier transform of $u$. One can also extend this operator to a wider space of functions
$$
{\cal L}_{\alpha}=\{u\left|\right.\int_{\mathbb{R}^n}\frac{|u(x)|}{1+|x|^{n+\alpha}}dx<\infty\}
$$
by
$$
<(-\lap)^{\alpha/2}u, \phi>=\int_{\mathbb{R}^n}u(-\lap)^{\alpha/2}\phi dx,\quad \text{for all} ~\phi\in C_0^\infty(\mathbb{R}^n).
$$
It is easy to verify that for $u\in{\cal L}_\alpha\cap C_{loc}^{1,1}$, the integral on the right hand side of (\ref{pp}) is well defined. Here, we consider the fractional Laplacian $(-\lap)^{\alpha/2}$ in this setting.

A great deal of results about symmetry and nonexistence of the solutions to elliptic equations involving $\lap$ were appeared (see \cite{CC1}, \cite{CFL}, \cite{ZLL} and the references therein). However, there are much fewer articles studying similar properties of solutions to equations involving $(-\lap)^{\alpha/2}$. The main difficulty is caused by the non-locality of the fractional Laplacian. To circumvent it, Caffarelli and Silvestre \cite{CS} introduced the extension method which reduced this nonlocal problem into a local one in higher dimensions. Chen, Li and Ou \cite{ChLO} provided the method of moving planes in integral forms, and treated the partial differential equations involving the fractional Laplacian by studying the equivalent integral equations. With these methods, a series of fruitful results have been obtained (see \cite{BCPS}, \cite{CFY}, \cite{ChL1}, \cite{ChL2}, \cite{ChZ}, \cite{FC}, \cite{ZC} and the references therein). For more articles concerning nonlocal equations, please see \cite{CDL}, \cite{CM}, \cite{FLS}, \cite{FW}, \cite{JW}, \cite{LQ}, \cite{Z}, \cite{ZCCY}, \cite{ZCY}, \cite{ZLCC} and the references therein.

Either by the \emph{extension} method or by integral equations, one may need to impose extra conditions on the solutions. But these extra conditions will not be necessary if we consider the pseudo differential equations, including equations involving $(-\lap)^{\alpha/2}$, directly.

Recently, a direct approach to carry on the method of moving planes for nonlocal problems on bounded or unbounded domains was developed by Chen, Li, and Li \cite{CLL}, which enables ones to systematically investigate equations involving the fractional Laplacian. In this approach, the \textit{narrow region principle} and the \textit{decay at infinity} for anti-symmetric functions will be repeatedly used. For readers' convenience, we will collect them in Section 2. Based on the approach, Chen, Li, and Li \cite{CLL} proved that if $u\in{\cal L}_\alpha\cap C_{loc}^{1,1}$ is a nonnegative solution of
\begin{equation}
(-\lap)^{\alpha/2}u=u^p(x),~~x\in \mathbb{R}^n,
\label{bb1}
\end{equation}
then

$(i)$ in the critical case $p=\frac{n+\alpha}{n-\alpha}$, $u$ is radially symmetric and monotone decreasing about some point;

$(ii)$ in the subcritical case $1<p<\frac{n+\alpha}{n-\alpha}$, $u\equiv0.\qquad\qquad\qquad\qquad\qquad\qquad\qquad\qquad\qquad\qquad\qquad\qquad$
To the Dirichlet problem on $\mathbb{R}^n_+$
\begin{equation}
\left\{\begin{array}{ll}
(-\lap)^{\alpha/2} u(x)=u^p(x), & \qquad x\in\mathbb{R}^n_+, \\
u(x)\equiv0, & \qquad x\notin\mathbb{R}^n_+,
\end{array}\right.
\label{bb2}\end{equation}
the authors of \cite{CLL} derived that if $1<p\leq\frac{n+\alpha}{n-\alpha}$ and $u\in{\cal L}_\alpha\cap C_{loc}^{1,1}$ is a nonnegative solution of (\ref{bb2}), then $u\equiv0$.

In this paper, we deal with the equation involving $(-\lap)^{\alpha/2}$ with more general nonlinearity. Throughout this paper, we suppose $u\in{\cal L}_\alpha\cap C_{loc}^{1,1}$.

Let us state the main results.

\begin{thm}
Suppose that $f(t)>0$ is strictly increasing in $t>0$, and $f(t)/t^{\frac{n+\alpha}{n-\alpha}}$ is non-increasing. Let $u$ be a positive solution of
\begin{equation}
(-\lap)^{\alpha/2} u(x)=f(u),\qquad x\in\mathbb{R}^n,
\label{01}
\end{equation}
then
\begin{equation}
u(x)=\frac{k}{(|x-x^o|^2+h^2)^{(n-\alpha)/2}}
\label{zlzlyy101}
\end{equation}
for some constants $k, h>0$ and some points $x^o\in\mathbb{R}^n$, and $f(t)=ct^{\frac{n+\alpha}{n-\alpha}}$ with $c>0$ for any $t\in(0, max_{\mathbb{R}^n}u]$.
\label{mthm1}
\end{thm}

\begin{remark}
(i) We know from Theorem \ref{mthm1} that $u$ is bounded on $\mathbb{R}^n$, and actually, only on $(0,max_{\mathbb{R}^n}]$ can we derive that $f(t)=ct^{\frac{n+\alpha}{n-\alpha}}$ with some $c>0$.

(ii) It sees in Theorem \ref{mthm1} that if $f(t)/t^{\frac{n+\alpha}{n-\alpha}}$ is strictly decreasing, then it contradict the conclusion that $f(t)=ct^{\frac{n+\alpha}{n-\alpha}}$ and hence (\ref{01}) possesses no positive solution.
\end{remark}

\begin{thm}
Suppose that $f(t)>0$ is strictly increasing in $t>0$ and $f(t)/t^{\frac{n+\alpha}{n-\alpha}}$ is non-increasing. Then
\begin{equation}
\left\{\begin{array}{ll}
(-\lap)^{\alpha/2} u(x)=f(u), & \qquad x\in\mathbb{R}^n_+, \\
u(x)\equiv0, & \qquad x\notin\mathbb{R}^n_+,
\end{array}\right.
\label{03}\end{equation}
possesses no positive solution.
\label{mthm2}
\end{thm}

To prove Theorems \ref{mthm1} and \ref{mthm2}, some new ideas are involved. Here we take the proof of Theorem \ref{mthm1} as an example to illustrate them. In the proof of Theorem \ref{mthm1}, to apply the method of moving planes, we need to make a proper Kelvin transform centered at $x^0\in\mathbb{R}^n$ ($x^0$ is arbitrarily chosen):
$$
v(x)=\frac{1}{|x-x^0|^{n-\alpha}}u(x^0+\frac{x-x^0}{|x-x^0|^2}),~~x\in\mathbb{R}^n\setminus\{x^0\},
$$
then $v(x)$ satisfies
\begin{equation}
(-\lap)^{\alpha/2} v(x)=\frac{1}{|x-x^0|^{n+\alpha}}f(|x-x^0|^{n-\alpha}v(x)),\quad x\in\mathbb{R}^n\setminus\{x^0\}.
\label{iyy3002}
\end{equation}
Through the moving plane method, we derive two possible cases: (i) $v$ is symmetric about some point $Q\neq x^0$, (ii) $v$ is symmetric about $x^0$. Here the new idea we want to underline is that in case (i), using the known fact that $u$ is bounded on $\mathbb{R}^n$, we derive by (\ref{iyy3002}) and the symmetry of $v$ that $f(t)=ct^{\frac{n+\alpha}{n-\alpha}}$ on $(0,max_{\mathbb{R}^n}]$ for some $c>0$, and reduce (\ref{01}) into
$$
(-\lap)^{\alpha/2}u(x)=Cu^{\frac{n+\alpha}{n-\alpha}}(x),~~~x\in\mathbb{R}^n,
$$
which leads to (\ref{zlzlyy101}) by the well known result \cite{ChLO}. To (ii), it immediately follows that there is no positive solution.

\begin{remark}
When $f(t)\equiv0$, we already derived the Liouville Theorem for $\alpha$-harmonic functions on $\mathbb{R}^n_+$ in our previous paper \cite{ZLCC}. The same result on $\mathbb{R}^n$ was established by Bogdan, Kulczycki and Nowak \cite{BKN}.
\end{remark}

\begin{remark}
$(i)$ Theorems \ref{mthm1} and \ref{mthm2} are generalizations of the elegant results in Bianchi \cite{BG} and Damascelli and Gladiali \cite{DG} from equations involving Laplacian to the corresponding equations involving the fractional Laplacian. It is worth noting that we drop the condition that $f$ is locally Lipschitz continuous. More specifically, we do not need the continuity of $f$ at all.

$(ii)$ In addition, we extend the results in \cite{CLL} to the problems with more general nonlinearity $f(u)$ on $\mathbb{R}^n$ and $\mathbb{R}^n_+$, and furthermore derive the exact representation of solution $u$ on $\mathbb{R}^n$.
\end{remark}

The organization of the paper is as follows. In Section 2, we list two important maximum principles for anti-symmetric functions which were given in \cite{CLL}. In Section 3, we prove Theorem \ref{mthm1}. Section 4 is devoted to the proof of Theorem \ref{mthm2}. In Section 5, we prove two claims which are used in Section 3.

\section{Two Known Maximum Principles}

Let $T$ be a hyperplane in $\mathbb{R}^n$ with the form
$$
T=\{x=(x_1,x')\in\mathbb{R}^n|x_1=\lambda ~for ~some~ \lambda\in\mathbb{R}\},
$$
where $x'=(x_2,x_3,\cdots,x_n)$. Denote the reflection of $x$ about the plane $T$ by
$$
\tilde{x}=(2\lambda-x_1, x_2, \cdots, x_n)
$$
and let
$$
H=\{x\in\mathbb{R}^n|x_1<\lambda\}~~and~~\tilde{H}=\{x|\tilde{x}\in H\}.
$$

We describe two key ingredients for the direct method of moving planes for the fractional Laplacian introduced in \cite{CLL}.

\begin{thm}(\cite{CLL},~Narrow Region Principle) Let $\Omega\subset H$ be a bounded narrow region contained in $\{x|\lambda-l<x_1<\lambda\}$ with small $l>0$. Suppose that $c(x)$ is bounded from below in $\Omega$, $u\in{\cal L}_\alpha\cap C_{loc}^{1,1}(\Omega)$ is lower semi-continuous on $\bar{\Omega}$ and satisfying
\begin{equation}
\left\{\begin{array}{lll}
(-\lap)^{\alpha/2}u(x)+c(x)u(x)\geq0, & \qquad in ~\Omega, \\
u(x)\geq0, & \qquad in ~H\backslash\Omega,\\
u(\tilde{x})=-u(x), & \qquad in ~H,
\end{array}\right.
\label{thm2.13}
\end{equation}
then we have
\begin{description}
\item[~~~\emph{(i)}]
\begin{equation}
u(x)\geq0~in~\Omega;
\label{thm2.14}
\end{equation}
\item[~~~\emph{(ii)}] furthermore, if $u=0$ at some point in $\Omega$, then
\begin{equation}
u(x)=0~~\text{almost everywhere in}~~\mathbb{R}^n;
\label{thm2.15555555555}
\end{equation}
\item[~~~\emph{(iii)}] conclusions (i) and (ii) hold for the unbounded region $\Omega$ if we further assume that
$$
\displaystyle\underset{|x| \ra \infty}{\underline{\lim}}u(x)\geq0.
$$
\end{description}
\label{thm2.2}
\end{thm}

\begin{thm}(\cite{CLL},~Decay at Infinity) Let $\Omega$ be an unbounded region in $H$. Assume that $u\in\mathcal{L}_\alpha\cap C_{loc}^{1,1}(\Omega)$ is a solution of
\begin{equation}
\left\{\begin{array}{lll}
(-\lap)^{\alpha/2}u(x)+c(x)u(x)\geq0, & \qquad in ~\Omega, \\
u(x)\geq0, & \qquad in ~H\backslash\Omega,\\
u(\tilde{x})=-u(x), & \qquad in ~H,
\end{array}\right.
\label{thm2.15}
\end{equation}
with
\begin{equation}
\displaystyle\underset{|x|\rightarrow\infty}{\underline{\lim}}|x|^\alpha c(x)\geq0,
\label{thm2.16}
\end{equation}
then there exists a constant $R_0>0$ (depending on $c(x)$ but independent of $u$) such that if
\begin{equation}
u(x^0)=\min_\Omega u(x)<0,
\label{thm2.17}
\end{equation}
then
$$
|x^0|\leq R_0.
$$
\label{thm2.3}
\end{thm}

\begin{remark} In Theorems \ref{thm2.2} and \ref{thm2.3}, the inequality
$$
(-\lap)^{\alpha/2}u(x)+c(x)u(x)\geq0,~~~\text{in} ~\Omega,
$$
and the condition (\ref{thm2.16}) are only required at points where $u$ is negative.
\label{rem2.1}
\end{remark}

\section{The proof of Theorem \ref{mthm1}}

In this section, we prove Theorem \ref{mthm1}.

\textbf{Proof of Theorem \ref{mthm1}.} Since there is no decay condition on $u$ near infinity, we are not able to carry the method of moving planes on $u$ directly. To circumvent this difficulty, we make a Kelvin transform. For any point $x^0\in\mathbb{R}^n$, let
$$
v(x)=\frac{1}{|x-x^0|^{n-\alpha}}u(x^0+\frac{x-x^0}{|x-x^0|^2}),~~x\in\mathbb{R}^n\setminus\{x^0\},
$$
be the Kelvin transform of $u$ centered at $x^0$. Then
\begin{equation}
(-\lap)^{\alpha/2} v(x)=\frac{1}{|x-x^0|^{n+\alpha}}f(|x-x^0|^{n-\alpha}v(x)),\quad x\in\mathbb{R}^n\setminus\{x^0\}.
\label{3002}
\end{equation}
The function $v$ is positive, decays to $0$ at infinity as $|x-x^0|^{\alpha-n}$ and may have a singularity at $x^0$. Next we prove that $v$ is radially symmetric, which implies the desired property of $u$.

For $x=(x_1,x_2,\cdots,x_n)$, choose any direction to be the $x_1$ direction. To prove that $v$ is radially symmetric, it suffices to show that $v$ is symmetric in $x_1$.

Write $x^0=(x^0_1,x^0_2,\cdots,x^0_n)$ and let that for $\lambda\in\mathbb{R}$ and $\lambda<x^0_1$,
$$
T_\lambda=\{x\in\mathbb{R}^n|x_1=\lambda\},~~x^\lambda=(2\lambda-x_1,x_2,\cdots,x_n),
$$
$$
\Sigma_\lambda=\{x\in\mathbb{R}^n|x_1<\lambda\},~~
v_\lambda(x):=v(x^\lambda),$$
and
$$
w_\lambda(x)=v_\lambda(x)-v(x),
~~\Sigma_\lambda^-=\{x\in\Sigma_\lambda|w_\lambda(x)<0\}.
$$
Then $w_\lambda(x)$ satisfies
\begin{eqnarray}
&&(-\lap)^{\alpha/2} w_\lambda(x)\nonumber\\&=&\frac{f(|x^\lambda-x^0|^{n-\alpha}v_\lambda(x))}{|x^\lambda-x^0|^{n+\alpha}}-
\frac{f(|x-x^0|^{n-\alpha}v(x))}{|x-x^0|^{n+\alpha}},\quad x\in\Sigma_\lambda\setminus\{(x^0)^\lambda\},
\label{32}
\end{eqnarray}
Noting
$$
\lim_{|x|\rightarrow\infty}w_\lambda(x)=0,
$$
we have that if $w_\lambda$ is negative somewhere in $\Sigma_\lambda$, then the negative minima of $w_\lambda$ are attained in the interior of $\Sigma_\lambda$. Since $f(t)>0$ for $t>0$ and $f(t)/t^{\frac{n+\alpha}{n-\alpha}}$ is non-increasing, it follows from (\ref{32}) that at points $x\in\Sigma_\lambda^-\setminus B_\epsilon((x^0)^\lambda)$ ($\epsilon$ small),
\begin{eqnarray}
&&(-\lap)^{\alpha/2}w_\lambda(x)\nonumber\\ &=&\frac{f(|x^\lambda-x^0|^{n-\alpha}v_\lambda(x))}{|x^\lambda-x^0|^{n+\alpha}}-
\frac{f(|x-x^0|^{n-\alpha}v(x))}{|x-x^0|^{n+\alpha}}\nonumber\\
&=&\frac{f(|x^\lambda-x^0|^{n-\alpha}v_\lambda(x))}{(|x^\lambda-x^0|^{n-\alpha}v_\lambda(x))
^{\frac{n+\alpha}{n-\alpha}}}v_\lambda^{\frac{n+\alpha}{n-\alpha}}(x)-
\frac{f(|x-x^0|^{n-\alpha}v(x))}{(|x-x^0|^{n-\alpha}v(x))^{\frac{n+\alpha}{n-\alpha}}}v^{\frac{n+
\alpha}{n-\alpha}}(x)\nonumber\\
&\geq&\frac{f(|x-x^0|^{n-\alpha}v(x))}{(|x-x^0|^{n-\alpha}v(x))
^{\frac{n+\alpha}{n-\alpha}}}\left(v_\lambda^{\frac{n+\alpha}{n-\alpha}}(x)-v^{\frac{n+
\alpha}{n-\alpha}}(x)\right)\nonumber\\
&=&\frac{n+\alpha}{n-\alpha}\cdot\frac{f(|x-x^0|^{n-\alpha}v(x))}{(|x-x^0|^{n-\alpha}v(x))
^{\frac{n+\alpha}{n-\alpha}}}\phi^{\frac{2\alpha}{n-\alpha}}(x)w_\lambda(x)\nonumber\\
&\geq&C\frac{f(|x-x^0|^{n-\alpha}v(x))}{(|x-x^0|^{n-\alpha}v(x))
^{\frac{n+\alpha}{n-\alpha}}}v^{\frac{2\alpha}{n-\alpha}}(x)w_\lambda(x),
\label{33}
\end{eqnarray}
where $v_\lambda(x)<\phi(x)<v(x)$. Here and below $C$ denotes a positive constant whose value may be different at different lines.

Denote $t=|x-x^0|^{n-\alpha}v(x)$, we know by the definition of $v(x)$ that $t\rightarrow u(x^0)>0$ as $|x|\rightarrow\infty$. Hence for any $x\in\Sigma_\lambda^-\setminus B_\epsilon((x^0)^\lambda)$, there exists a constant $c_0>0$ such that $t\geq c_0$. Using that $f(t)/t^{\frac{n+\alpha}{n-\alpha}}$ is non-increasing, it yields
$$\frac{f(t)}{t^{\frac{n+\alpha}{n-\alpha}}}\leq\frac{f(c_0)}{c_0^{\frac{n+\alpha}{n-\alpha}}},$$
and by (\ref{33}),
\begin{eqnarray}
(-\lap)^{\alpha/2}w_\lambda(x)
&\geq&C\frac{f(|x-x^0|^{n-\alpha}v(x))}{(|x-x^0|^{n-\alpha}v(x))
^{\frac{n+\alpha}{n-\alpha}}}v^{\frac{2\alpha}{n-\alpha}}(x)w_\lambda(x)\label{j33}\\
&\geq&Cv^{\frac{2\alpha}{n-\alpha}}(x)w_\lambda(x).
\end{eqnarray}
Therefore,
\begin{equation}
(-\lap)^{\alpha/2}w_\lambda(x)+c(x)w_\lambda(x)\geq0,\quad x\in\Sigma_\lambda^-\setminus B_\epsilon((x^0)^\lambda)
\label{34}
\end{equation}
with
\begin{equation}
c(x)=-Cv^{\frac{2\alpha}{n-\alpha}}(x).
\label{35}
\end{equation}

\textbf{Step 1.} We prove that for $\lambda<0$, $|\lambda|$ large enough,
\begin{equation}
w_\lambda\geq0~~ \text{in}~~ \Sigma_\lambda\setminus\{(x^0)^\lambda\}
\label{y1}
\end{equation}
by using Theorem \ref{thm2.3}.

Let us first admit:

\textbf{Claim 3.1.} \textsl{For $\lambda$ sufficiently negative, there exists $\epsilon>0$ and $c_\lambda>0$, such that}
\begin{equation}
w_\lambda(x)\geq c_\lambda,~~~~x\in B_\epsilon((x^0)^\lambda)\setminus\{(x^0)^\lambda\}.
\label{0360}
\end{equation}

Its proof will be given in Section 5. From (\ref{0360}), one can see that $\Sigma_\lambda^-$ has no intersection with $B_\epsilon((x^0)^\lambda)$. Using (\ref{35}), it is easy to verify that for $|x|$ sufficiently large,
\begin{equation}
c(x)\sim\frac{1}{|x|^{2\alpha}},
\label{36}
\end{equation}
hence $c(x)$ satisfies (\ref{thm2.16}). Applying Theorem \ref{thm2.3} to $w_\lambda$ with
$$
H=\Sigma_\lambda~~\text{and}~~\Omega=\Sigma_\lambda^-,
$$
we conclude that there exists $R_0>0$ (independent of $\lambda$), such that if $\bar{x}$ is a negative minimum of $w_\lambda$ in $\Sigma_\lambda$, then
\begin{equation}
|\bar{x}|\leq R_0.
\label{37}
\end{equation}
Now for $\lambda\leq-R_0$, we have (\ref{y1}).

\textbf{Step 2.} \emph{Step 1} provides a starting point such that we can now move the plane $T_\lambda$ to the right as long as (\ref{y1}) holds to its limiting position. Denoting
$$
\lambda_0=\sup\{\lambda<x^0_1|w_\mu(x)\geq0, ~x\in\Sigma_\mu\setminus\{(x^0)^\mu\},~\mu\leq\lambda\},
$$
we show
\begin{equation}
w_{\lambda_0}(x)\equiv0, ~~~x\in \Sigma_{\lambda_0}\setminus\{(x^0)^{\lambda_0}\}.
\label{39}
\end{equation}
To do so, let us consider two possibilities: (i)$~\lambda_0<x^0_1$; and (ii)$~\lambda_0=x^0_1$.

\textbf{(\textbf{i}).} For $\lambda_0<x^0_1,$ we deduce that if
\begin{equation}
w_{\lambda_0}(x)\not\equiv0,~~x\in \Sigma_{\lambda_0}\setminus\{(x^0)^{\lambda_0}\},
\label{z1}
\end{equation}
then the plane $T_\lambda$ can be moved further to the right, to be more rigorous, there exists some $\varepsilon>0$ such that for any $\lambda\in(\lambda_0,\lambda_0+\varepsilon)$,
\begin{equation}
w_{\lambda}(x)\geq0, ~~~x\in \Sigma_{\lambda}\setminus\{(x^0)^\lambda\}.
\label{310}
\end{equation}
This is a contradiction with the definition of $\lambda_0$. Hence we must have (\ref{39}).

\textbf{(\textbf{i})$_1$}. Now we prove (\ref{310}) by combining Theorems \ref{thm2.2} and \ref{thm2.3}. Let us admit:

\textbf{Claim 3.2.} \textsl{There exists $c_0>0$ such that for sufficiently small $\eta$,
\begin{equation}
w_{\lambda_0}(x)\geq c_0,~~x\in B_\eta((x^0)^{\lambda_0})\setminus\{(x^0)^{\lambda_0}\}.
\label{y2}
\end{equation}}

We put its proof in Section 5. From (\ref{y2}), one knows that $\Sigma_{\lambda_0}^-$ does not intersect with $B_\eta((x^0)^{\lambda_0})$, hence for $\lambda$ sufficiently close to $\lambda_0$, $\Sigma_{\lambda}^-$ has no intersection with $B_\eta((x^0)^{\lambda})$ either. To see this, for any $x\in B_\eta((x^0)^{\lambda})\setminus\{(x^0)^{\lambda}\}$, we find $\bar{x}\in B_\eta((x^0)^{\lambda_0})\setminus\{(x^0)^{\lambda_0}\}$ such that $x^\lambda=\bar{x}^{\lambda_0}$. Then, for $\lambda$ sufficiently close to $\lambda_0$, since $v(x)$ is continuous in $\Sigma_\lambda$, it follows
\begin{eqnarray}
w_\lambda(x)&=&[v_\lambda(x)-v_{\lambda_0}(\bar{x})]+[v_{\lambda_0}(\bar{x})-v(\bar{x})]+[v(\bar{x})-v(x)]\nonumber\\
&\geq&c_0+[v(\bar{x})-v(x)]\nonumber\\
&\geq&0,~~x\in B_\eta((x^0)^{\lambda})\setminus\{(x^0)^{\lambda}\},
\label{y3}
\end{eqnarray}
which verifies that $\Sigma_{\lambda}^-$ does not intersect with $B_\eta((x^0)^{\lambda})$.

By (\ref{37}), the negative minimum of $w_\lambda$ connot be attained outside of $B_{R_0}(0)$. Next we will prove that it can neither be attained inside of $B_{R_0}(0)$, i.e., for $\lambda$ sufficiently close to $\lambda_0$,
\begin{equation}
w_\lambda(x)\geq0,~~x\in(\Sigma_\lambda\cap B_{R_0}(0))\setminus\{(x^0)^\lambda\}.
\label{312}
\end{equation}

Actually, using Theorem \ref{thm2.2}, there is a small $\delta>0$, such that for $\lambda\in[\lambda_0,\lambda_0+\delta)$, if
\begin{equation}
w_\lambda(x)\geq0,~~x\in\Sigma_{\lambda_0-\delta}\setminus\{(x^0)^{\lambda}\},
\label{313}
\end{equation}
then
\begin{equation}
w_\lambda(x)\geq0,~~x\in(\Sigma_\lambda\setminus\Sigma_{\lambda_0-\delta})\setminus\{(x^0)^{\lambda}\}.
\label{314}
\end{equation}
To derive (\ref{314}) from (\ref{313}), we use Theorem \ref{thm2.2} with
$$
H=\Sigma_\lambda~\text{and the narrow region}~\Omega=(\Sigma_\lambda^-\setminus\Sigma_{\lambda_0-\delta}),
$$
while the lower bound of $c(x)$ can be seen from (\ref{35}) and (\ref{36}).

Now if (\ref{313}) holds, then one can immediately get (\ref{312}) from (\ref{313}) and (\ref{314}), so what left is to show (\ref{313}), and by (\ref{37}) and (\ref{y3}) we only need to prove
\begin{equation}
w_\lambda(x)\geq0,~~x\in(\Sigma_{\lambda_0-\delta}\cap B_{R_0}(0))\setminus B_\eta((x^0)^\lambda).
\label{315}
\end{equation}
In fact, when $\lambda_0<x_1^0$, we have
\begin{equation}
w_{\lambda_0}(x)>0,~~x\in\Sigma_{\lambda_0}\setminus\{(x^0)^{\lambda_0}\}.
\label{316}
\end{equation}
If not, then there exists some $\hat{x}$ such that
$$
w_{\lambda_0}(\hat{x})=\min_{\Sigma_{\lambda_0}}w_{\lambda_0}(x)=0
$$
and it implies
\begin{eqnarray}
&&(-\lap)^{\alpha/2}w_\lambda(\hat{x})\nonumber\\
&=&C_{n,\alpha}PV\int_{\mathbb{R}^n}\frac{-w_{\lambda_0}(y)}{|\hat{x}-y|^{n+\alpha}}dy\nonumber\\
&=&C_{n,\alpha}PV\int_{\Sigma_{\lambda_0}}\frac{-w_{\lambda_0}(y)}{|\hat{x}-y|^{n+\alpha}}dy
+\int_{\mathbb{R}^n\setminus\Sigma_{\lambda_0}}\frac{-w_{\lambda_0}(y)}{|\hat{x}-y|^{n+\alpha}}dy
\nonumber\\
&=&C_{n,\alpha}PV\int_{\Sigma_{\lambda_0}}\frac{-w_{\lambda_0}(y)}{|\hat{x}-y|^{n+\alpha}}dy
+\int_{\Sigma_{\lambda_0}}\frac{w_{\lambda_0}(y)}{|\hat{x}-y^\lambda|^{n+\alpha}}dy
\nonumber\\
&=&C_{n,\alpha}PV\int_{\Sigma_{\lambda_0}}\left(\frac{1}{|\hat{x}-y^\lambda|^{n+\alpha}}-
\frac{1}{|\hat{x}-y|^{n+\alpha}}\right)w_{\lambda_0}(y)dy\nonumber\\
&\leq&0.
\label{317}
\end{eqnarray}
On the other hand, we have by the monotonicity of $f(t)/t^{\frac{n+\alpha}{n-\alpha}}$ that
\begin{eqnarray}
&&(-\lap)^{\alpha/2}w_\lambda(\hat{x})\nonumber\\ &=&\frac{f(|\hat{x}^\lambda-x^0|^{n-\alpha}v_\lambda(\hat{x}))}{|\hat{x}^\lambda-x^0|^{n+\alpha}}-
\frac{f(|\hat{x}-x^0|^{n-\alpha}v(\hat{x}))}{|\hat{x}-x^0|^{n+\alpha}}\nonumber\\
&=&\frac{f(|\hat{x}^\lambda-x^0|^{n-\alpha}v_\lambda(\hat{x}))}{(|\hat{x}^\lambda-x^0|^{n-\alpha}v_\lambda(\hat{x}))
^{\frac{n+\alpha}{n-\alpha}}}v_\lambda^{\frac{n+\alpha}{n-\alpha}}(\hat{x})-
\frac{f(|\hat{x}-x^0|^{n-\alpha}v(\hat{x}))}{(|\hat{x}-x^0|^{n-\alpha}v(\hat{x}))^{\frac{n+\alpha}{n-\alpha}}}v^{\frac{n+
\alpha}{n-\alpha}}(\hat{x})\nonumber\\
&=&\left(\frac{f(|\hat{x}^\lambda-x^0|^{n-\alpha}v(\hat{x}))}{(|\hat{x}^\lambda-x^0|^{n-\alpha}v(\hat{x}))
^{\frac{n+\alpha}{n-\alpha}}}-
\frac{f(|\hat{x}-x^0|^{n-\alpha}v(\hat{x}))}{(|\hat{x}-x^0|^{n-\alpha}v(\hat{x}))^{\frac{n+\alpha}{n-\alpha}}}\right)v^{\frac{n+
\alpha}{n-\alpha}}(\hat{x})\nonumber\\
&\geq&0.
\label{318}
\end{eqnarray}
From (\ref{317}) and (\ref{318}), one derives that $w_{\lambda_0}(x)\equiv0~\text{in}~\Sigma_{\lambda_0}\setminus\{(x^0)^{\lambda_0}\}$, which is a contradiction with (\ref{z1}). This proves (\ref{316}).

It follows from (\ref{316}) that there exists a constant $c_0>0$, such that
$$
w_{\lambda_0}(x)\geq c_0,~~x\in\overline{(\Sigma_{\lambda_0-\delta}\cap B_{R_0}(0))\setminus B_\eta((x^0)^{\lambda_0})}.
$$
Since $w_\lambda$ depends continuously on $\lambda$, there exists $\varepsilon>0$ and $\varepsilon<\delta$, such that for all $\lambda\in(\lambda_0,\lambda_0+\varepsilon)$, we have
\begin{equation}
w_{\lambda}(x)\geq0,~~x\in\overline{(\Sigma_{\lambda_0-\delta}\cap B_{R_0}(0))\setminus B_\eta((x^0)^{\lambda})}.
\label{319}
\end{equation}
This verifies (\ref{315}) and hence (\ref{313}) is also proved.

Combining (\ref{37}), (\ref{y3}), (\ref{314}) and (\ref{319}), we conclude (\ref{310}) for all $\lambda\in(\lambda_0,\lambda_0+\varepsilon)$, which is a contradiction with the definition of $\lambda_0$. Therefore, we must have (\ref{39}).

\textbf{(i)$_2$}. Now we know by (\ref{39}) that $v(x)$ is radially symmetric about some point $Q$ different from $x^0$ and $v(x)\not\equiv c$ is bounded near $x^0$, thus $u(x)\not\equiv c$ is bounded in $\mathbb{R}^n$ too. And also, $(-\lap)^{\alpha/2}v(x)$ is radially symmetric about the same point $Q$, which can be easily proved through elementary computation with the help of (\ref{pp}). Thus the right hand side of the following equation
\begin{equation}
(-\lap)^{\alpha/2} v(x)=\frac{1}{|x-x^0|^{n+\alpha}}f(|x-x^0|^{n-\alpha}v(x)),\quad x\in\mathbb{R}^n\setminus\{x^0\}
\label{ref1}
\end{equation}
should have the same symmetry. From this, we are able to prove that
\begin{equation}
f(|x-x^0|^{n-\alpha}v(x))=C\left[|x-x^0|^{n-\alpha}v(x)\right]^{\frac{n+\alpha}{n-\alpha}},\quad x\in\mathbb{R}^n\setminus\{x^0\},
\label{eref1}
\end{equation}
for some positive constant $C$ ($C$ cannot be $0$, otherwise $v(x)\equiv c>0$ is radially symmetric about $x^0$, refer to \cite{CDL}). Yet the proof of (\ref{eref1}) is sophisticated, and we postpone the proof for a while.

Now set
$$t(x)=|x-x^0|^{n-\alpha}v(x),\quad x\in\mathbb{R}^n\setminus\{x^0\}.$$
Since $u$ is bounded on $\mathbb{R}^n$, it infers by
$$|x-x^0|^{n-\alpha}v(x)=u(x^0+\frac{x-x^0}{|x-x^0|^2})$$
that $t\in(0, max_{\mathbb{R}^n}u]$. Noting (\ref{eref1}), we get $f(t)=Ct^{\frac{n+\alpha}{n-\alpha}}$ for any $t\in(0, max_{\mathbb{R}^n}u]$ and some positive constant $C$ and (\ref{01}) becomes
$$
(-\lap)^{\alpha/2}u(x)=Cu^{\frac{n+\alpha}{n-\alpha}}(x),~~~x\in\mathbb{R}^n,
$$
which gives by \cite{ChLO} that
\begin{equation}
u(x)=\frac{k}{(|x-x^o|^2+h^2)^{(n-\alpha)/2}}
\end{equation}
for some constants $k, h>0$ and some point $x^o\in\mathbb{R}^n$.

(\textbf{i})$_3$. Now let us prove (\ref{eref1}). We turn (\ref{ref1}) into the following form
\begin{equation}
[(-\lap)^{\alpha/2}v(x)]v^{\frac{n-\alpha}{n+\alpha}}=\frac{f(|x-x^0|^{n-\alpha}v(x))}{[|x-x^0|^{n-\alpha}v(x)]^{\frac{n+\alpha}{n-\alpha}}},\quad x\in\mathbb{R}^n\setminus\{x^0\}.
\label{ref2}
\end{equation}
Since the left hand side of (\ref{ref2}) is radially symmetric about $Q$, it implies that
$$
\frac{f(|x-x^0|^{n-\alpha}v(x))}{[|x-x^0|^{n-\alpha}v(x)]^{\frac{n+\alpha}{n-\alpha}}}=g(|x-Q|),\quad x\in\mathbb{R}^n\setminus\{x^0\},
$$
is also radially symmetric about $Q$ for some proper function $g$. For simplicity, we denote
$$
F(t):=\frac{f(t)}{t^{\frac{n+\alpha}{n-\alpha}}},\qquad t=t(x)=|x-x^0|^{n-\alpha}v(x).
$$

Let $l_{Qx^0}$ be the line passing through points $Q$ and $x^0$, and $B_R(Q)$ the sphere centered at $Q$ with radius $R$, then for any $x\in B_R(Q)$, one can easily see that the maximum point and the minimum point of $t(x)$ are located on the line $l_{Qx^0}$. Suppose that the maximum point of $t(x)$ on $\mathbb{R}^n$ is $\tilde{P}$, then $\tilde{P}\in l_{Qx^0}$ ($\tilde{P}\neq Q$ and $\tilde{P}\neq x^0$), and $\tilde{P}$, $x^0$ are at the different sides of $Q$ on $l_{Qx^0}$.

Define $\tilde{x}:=2Q-x$, then the symmetry point of $\tilde{P}$ about $Q$ is $P:=2Q-\tilde{P}$. Let $s_{x^0P}$ be the line segment between points $x^0$ and $P$, $P\in s_{x^0P}$, but $x^0\not\in s_{x^0P}$. Then for any $A\in s_{x^0P}$, we have $t(A)<t(\tilde{A})$ and $t(\tilde{A})-t(A)\geq c_0$ for some positive constant $c_0$. Without loss of generality, we may assume $t(A)=m$ and $t(\tilde{A})=M$. Then, because $g(A)=g(\tilde{A})$, so $F(m)=F(M)$. Using that $F(t)$ is non-increasing in $t$ on $(0,+\infty)$, we derive
\begin{equation}
F(t)=C_0>0,~~~\text{on}~[m,M]~\text{with}~M-m\geq c_0>0.
\label{min0}
\end{equation}

Concretely, when $A\rightarrow x^0$, one has $t(A)\rightarrow0$ and thus $\tilde{A}\rightarrow\tilde{x^0}$, $t(\tilde{A})\rightarrow t(\tilde{x^0}):=a$, $a>0$. Then by the continuity of $t(x)$ and the above technique, we have
\begin{equation}
F(t)=C_1>0,~~~\text{on}~(0,a).
\label{min1}
\end{equation}
When $A=P$, suppose $t(P)=b$, we already know $t(\tilde{P})=max_{\mathbb{R}^n}u$. Then repeat the above process, we have
\begin{equation}
F(t)=C_2>0,~~~\text{on}~[b,max_{\mathbb{R}^n}u].
\label{min2}
\end{equation}

If $a\geq b$, then by (\ref{min1}), (\ref{min2}) and the continuity of $F(t)$ for $t>0$, we already proved that
\begin{equation}
F(t)=C>0,~~~\text{on}~(0,max_{\mathbb{R}^n}u].
\label{min3}
\end{equation}

If $a<b$, we can also show that (\ref{min3}) is still true. In fact, let us seek a finite sequence of points $\{A_i\}\subset s_{x^0P}$, $i=1,2,\cdots,k$, such that
$$
t(A_{i+1})=t(\tilde{A_i}):=m_i, ~i=1,2,\cdots,k-1,
$$
$$
t(A_1)=a:=m_0,~~t(\tilde{A_k}):=m_k~~\text{and}~~m_{k-1}\leq b\leq m_k.
$$
Since $f(t)\in(0,b]$ is continuous on $s_{x^0P}$, and $m_{i}-m_{i-1}\geq c_0$ which can be seen from (\ref{min0}), $i=1,2,\cdots,k$, we can find such $\{A_i\}$ indeed, and we also know that $k\leq\lfloor\frac{b-a}{c_0}\rceil+1$, where $\lfloor\frac{b-a}{c_0}\rceil$ denotes the integer part of $\frac{b-a}{c_0}$. It follows from (\ref{min0}) that
\begin{equation}
F(t)=c_i>0,~~~\text{on}~[m_i,m_{i+1}],~~i=0,1,2,\cdots,k.
\label{min4}
\end{equation}
According to
$$
(0,a)\cup[a,m_1]\cup[m_1,m_2]\cup\cdots\cup[m_{k-1},m_k]\cup[b,max_{\mathbb{R}^n}u]=(0,max_{\mathbb{R}^n}u],
$$
then by (\ref{min1}), (\ref{min2}), (\ref{min4}) and the continuity of $F(t)$ on $t>0$, we can also derive (\ref{min3}) for $a<b$.

Hence (\ref{eref1}) is proved.

Till now, we have proved Theorem \ref{mthm1} under Possibility $(i)$.

\textbf{(\textbf{ii}).} From the definition of $\lambda_0$, one knows that
\begin{equation}
\lambda_0=x_1^0~\text{and}~w_{\lambda_0}\geq0,~~x\in\Sigma_{\lambda_0}\setminus\{(x^0)^{\lambda_0}\}.
\label{327}
\end{equation}
Now we move the plane $T_\lambda$ from $+\infty$ to the left and similarly derive that either

\begin{equation}
\lambda_0>x_1^0 ~\text{and}~w_{\lambda_0}\equiv0,~~x\in\Sigma_{\lambda_0}.\qquad\quad\quad\qquad
\label{329}
\end{equation}
or
\begin{equation}
\lambda_0=x_1^0 ~\text{and}~w_{\lambda_0}\leq0,~~x\in\Sigma_{\lambda_0}.\qquad\quad\quad\qquad
\label{328}
\end{equation}

The case described by $(\ref{329})$ can be handled with the same way as (\textbf{i}). Now from (\ref{327}) and (\ref{328}), we have
$$
\lambda_0=x_1^0 ~\text{and}~w_{\lambda_0}\equiv0,~x\in\Sigma_{\lambda_0}.
$$

So far, we have proved that $v$ is symmetric about the plane $T_{x^0_1}$. Since the $x_1$ direction can be chosen arbitrarily, we have actually shown that $v$ is radially symmetric about $x^0$.

For any two points $X^i\in\mathbb{R}^n,~i=1,2$, choose $x^0=\frac{X^1+X^2}{2}$. Since $v$ is radially symmetric about $x^0$, so is $u$, hence $u(X^1)=u(X^2)$. This implies that $u$ is a constant, without loss of generality, we may suppose $u(x)=a_0>0$. Now (\ref{01}) holds if and only if $f(a_0)=0$, this is a contradiction with the condition ``$f(t)>0$ for $t>0$". Hence (\ref{01}) does not exist any positive solution.

The proof of Theorem \ref{mthm1} ends here.

\section{The Proof of Theorem \ref{mthm2}}

We first indicate two propositions.
\begin{proposition}
~~Assume that $u\in\mathcal{L}_\alpha\cap C_{loc}^{1,1}(\mathbb{R}^n_+)$ is a locally bounded positive solution to the problem
\begin{equation}
\left\{\begin{array}{ll}
(-\lap)^{\alpha/2} u(x)=f(u), & \qquad x\in\mathbb{R}^n_+, \\
u(x)\equiv0, & \qquad x\notin\mathbb{R}^n_+,
\end{array}\right.
\label{423}
\end{equation}
then it is also a solution to the integral equation
\begin{equation}
u(x)=\int_{\mathbb{R}^n_+}G_\infty(x,y)f(u(y))dy;
\label{424}
\end{equation}
and vice versa. Here $G_\infty(x,y)$ is the Green function to (\ref{423}):
$$
G_\infty(x,y)=\frac{A_{n,\alpha}}{s^{\frac{n-\alpha}{2}}}\left[1-B\frac{1}{(t+s)^{\frac{n-2}{2}}}
\int_0^{\frac{s}{t}}\frac{(s-tb)^{\frac{n-2}{2}}}{b^{\alpha/2}(1+b)}db\right],
$$
where $A_{n,\alpha}, B$ are positive constants, $t=4x_ny_n$ and $s=|x-y|^2$.
\label{pro4.123}
\end{proposition}

The proof of Proposition \ref{pro4.123} is similar to Theorem 4.1 in \cite{CFY}. Since the condition ``$f(t)>0$ is strictly increasing for $t>0$" guarantees that $f(t)$ is locally bounded on $(0,+\infty)$ and $f(t)\geq C_0$ for $t>R$, where $R>0$ is sufficiently large and $C_0$ is a positive constant, we only need to substitute $f(u)$ for $u^p$ in the proof of \cite{CFY}. We omit the details.

\begin{proposition}
(\cite{CFY})
~~If $\frac{t}{s}$ is sufficiently small, then for any $x=(x',x_n)$, $y=(y',y_n)\in \mathbb{R}^n_+$,
\begin{equation}
\frac{c_{n,\alpha}}{s^{(n-\alpha)/2}}\frac{t^{\alpha/2}}{s^{\alpha/2}}\leq G_\infty(x,y)\leq\frac{C_{n,\alpha}}{s^{(n-\alpha)/2}}\frac{t^{\alpha/2}}{s^{\alpha/2}},
\label{425}
\end{equation}
that is
\begin{equation}
G_\infty(x,y)\sim\frac{t^{\alpha/2}}{s^{n/2}},
\label{426}
\end{equation}
where $t=4x_ny_n,~s=|x-y|^2,$
$c_{n,\alpha}$ and $C_{n,\alpha}$ stand for different positive constants and only depend on $n$ and $\alpha$.
\label{pro4CFYG}
\end{proposition}

Now let us prove Theorem \ref{mthm2}.

\textbf{Proof of Theorem \ref{mthm2}.} Let us make a Kelvin transform again:
$$
v(x)=\frac{1}{|x-x^0|^{n-\alpha}}u(x^0+\frac{x-x^0}{|x-x^0|^2}),~~x\in\mathbb{R}^n_+,
$$
where the center $x^0$ is on the boundary $\partial\mathbb{R}^n_+$ in order that $\mathbb{R}^n_+$ is invariant under the transform. Obviously, $v$ satisfies
\begin{equation}
(-\lap)^{\alpha/2} v(x)=\frac{1}{|x-x^0|^{n+\alpha}}f(|x-x^0|^{n-\alpha}v(x)),\quad x\in\mathbb{R}^n_+,
\label{401}
\end{equation}
and is positive, decays to $0$ at infinity as $|x-x^0|^{\alpha-n}$ and may have a singularity at $x^0$. For $x=(x_1,x_2,\cdots,x_n)$, choose any direction in $\mathbb{R}^{n-1}$(the boundary of $\mathbb{R}^n_+$) to be $x_1$ direction. We first prove
\begin{equation}
v=v(x_n),
\label{44401}
\end{equation}
that is, $v$ is rotational symmetric about the $x_n$ axis. Next we derive a contradiction by (\ref{44401}) and prove that any positive solution of (\ref{03}) does not exist.

To prove (\ref{44401}), it suffices to show that $v$ is axisymmetric about any line parallel to $x_n$ axis. Write $x^0=(x^0_1,x^0_2,\cdots,x^0_n)$, for $\lambda<x^0_1$, $\lambda\in\mathbb{R}$, and for any $x\in\mathbb{R}^n_+$, denote
$$
T_\lambda=\{x\in\mathbb{R}^n|x_1=\lambda\},~~x^\lambda=(2\lambda-x_1,x_2,\cdots,x_n),
$$
$$
\Sigma_\lambda=\{x\in\mathbb{R}^n_+|x_1<\lambda\},~~
v_\lambda(x):=v(x^\lambda),
$$
and
$$
w_\lambda(x)=v_\lambda(x)-v(x),~~
\Sigma_\lambda^-=\{x\in\Sigma_\lambda|w_\lambda(x)<0\}.
$$
Then $w_\lambda(x)$ satisfies
\begin{equation}
(-\lap)^{\alpha/2} w_\lambda(x)=\frac{f(|x^\lambda-x^0|^{n-\alpha}v_\lambda(x))}{|x^\lambda-x^0|^{n+\alpha}}-
\frac{f(|x-x^0|^{n-\alpha}v(x))}{|x-x^0|^{n+\alpha}},~~x\in\Sigma_\lambda.
\label{402}
\end{equation}

Since $f(t)>0$ for $t>0$, and $f(t)/t^{\frac{n+\alpha}{n-\alpha}}$ is non-increasing, similarly to (\ref{34}) and (\ref{35}) in the previous section, we have by (\ref{402}) that for $\epsilon$ small,
\begin{equation}
(-\lap)^{\alpha/2}w_\lambda(x)+c(x)w_\lambda(x)\geq0,\quad x\in\Sigma_\lambda^-\setminus B_\epsilon((x^0)^\lambda),
\label{403}
\end{equation}
where
\begin{equation}
c(x)=-Cv^{\frac{2\alpha}{n-\alpha}}(x).
\label{404}
\end{equation}
Hence the lower bound of $c(x)$ in $\Sigma_\lambda^-\setminus B_\epsilon((x^0)^\lambda)$ can be obtained from (\ref{404}), and it follows by the asymptotic behavior
$$
v(x)\sim\frac{1}{|x|^{n-\alpha}},~~\text{for}~|x|~\text{large}
$$
that
$$
\lim_{|x|\rightarrow\infty}w_\lambda(x)=0~\text{and}~c(x)\sim\frac{1}{|x|^{2\alpha}}~\text{for}~|x|~\text{large}.
$$
These guarantee that we can apply Theorems \ref{thm2.2} and \ref{thm2.3} to show the following:

\emph{Step 1.} For $\lambda$ sufficiently negative,
$$
w_\lambda(x)\geq0, ~~x\in\Sigma_\lambda.
$$

\emph{Step 2.} Denoting
$$
\lambda_0=\sup\{\lambda|w_\mu(x)\geq0,x\in\Sigma_\mu;\mu\leq\lambda<x_1^0\},
$$
we have
\begin{equation}
w_{\lambda_0}(x)\equiv0,~~x\in\Sigma_{\lambda_0}.
\label{resul11}
\end{equation}

The proofs of two steps are in general similar to the proof of Theorem \ref{mthm1}, here we only state the differences.

To prove (\ref{resul11}), we also need to consider two possibilities: (i) $\lambda_0<x^0_1$; and (ii) $\lambda_0=x^0_1.$

\textbf{(\textbf{i}).} A similar proof obtaining (\ref{39}) leads to
\begin{equation}
\lambda_0<x^0_1,~~w_{\lambda_0}\equiv0,~~x\in\Sigma_{\lambda_0}.
\label{405}
\end{equation}
By $l_Q$ we denote the line parallel to $x_n$ axis passing through the point $Q\in\partial\mathbb{R}^n_+$. For any point $P=(p_1,\cdots,p_n)$, denote $p'=(p_1,\cdots,p_{n-1})$.

Since one can choose the $x_1$ direction arbitrarily, it easily derives by (\ref{405}) that $v(x)$ is axisymmetric about the line $l_Q$ different from $l_{x^0}$. Thus similarly to (\ref{ref1}), the right hand side of the following equation
$$
(-\lap)^{\alpha/2} v(x)=\frac{1}{|x-x^0|^{n+\alpha}}f(|x-x^0|^{n-\alpha}v(x)),\quad x\in\mathbb{R}^n_+,
$$
should also be axisymmetric about line $l_Q$. Since $l_Q\neq l_{x^0}$, similarly to (\ref{eref1}), we have
$$
f(|x-x^0|^{n-\alpha}v(x))=C\left[|x-x^0|^{n-\alpha}v(x)\right]^{\frac{n+\alpha}{n-\alpha}},
$$
for some positive constants $C$ (note $C\neq0$, otherwise $v(x)=v(x_n)$ is axisymmetric about $l_{x^0}$, refer to \cite{ZLCC}).

Till now, we have shown that $v(x)$ is bounded near $x^0$ and $f(t)=Ct^{\frac{n+\alpha}{n-\alpha}}$ for some positive constants $C$. Hence
\begin{equation}
u(x)\sim\frac{1}{|x|^{n-\alpha}}~~~\text{near infinity}
\label{ii408}
\end{equation}
and (\ref{03}) becomes
\begin{equation}
\left\{\begin{array}{ll}
(-\lap)^{\alpha/2} u(x)=Cu^{\frac{n+\alpha}{n-\alpha}}(x), & \qquad x\in\mathbb{R}^n_+, \\
u(x)\equiv0, & \qquad x\notin\mathbb{R}^n_+.
\end{array}\right.
\label{ii409}\end{equation}
The result of section 3.2.1 in \cite{CLL} yields that (\ref{03}) possesses no positive solution.

\textbf{(\textbf{ii}).} From the definition of $\lambda_0$, we know
\begin{equation}
\lambda_0=x_1^0~\text{and}~w_{\lambda_0}\geq0,~~x\in\Sigma_{\lambda_0}.
\label{327777}
\end{equation}
Now let us move the plane $T_\lambda$ from $+\infty$ to the left. Similarly, we can derive that either
\begin{equation}
\lambda_0>x_1^0 ~\text{and}~w_{\lambda_0}\equiv0,~~x\in\Sigma_{\lambda_0}
\label{407}
\end{equation}
or
\begin{equation}
\lambda_0=x_1^0 ~\text{and}~w_{\lambda_0}\leq0,~~x\in\Sigma_{\lambda_0}.
\label{408}
\end{equation}

The case (\ref{407}) can be handled with the same way as in (\textbf{i}). Now from (\ref{327777}) and (\ref{408}), we get
$$
\lambda_0=x_1^0 ~\text{and}~w_{\lambda_0}\equiv0,~~x\in\Sigma_{\lambda_0}.
$$
So far, we have proved that $v$ is symmetric about the plane $T_{x^0_1}$. Noting that the $x_1$ direction can be chosen arbitrarily, we have actually shown that $v$ is axially symmetric about the line parallel to $x_n$ axis and passing through $x^0$. Because $x^0$ is any point on $\partial\mathbb{R}^n_+$, we already deduce that the original solution $u$ is independent of the first $n-1$ variables, that is $u=u(x_n)$.

Next, we will show that this will contradict with the finiteness of $u$. By Proposition \ref{pro4.123}, we only need to show that $u=u(x_n)$ contradicts the finiteness of the integral
$$
\int_{\mathbb{R}^n_+}G_\infty(x,y)f(u(y))dy.
$$
In fact, by Proposition \ref{pro4CFYG}, if $u(x)=u(x_n)$ is a solution of (\ref{424}), then for each fixed $x\in\mathbb{R}^n_+$ and $R$ large enough, we have by (\ref{425}),
\begin{eqnarray}
+\infty>u(x)&=&\int_0^\infty f(u)\int_{\mathbb{R}^{n-1}}G_\infty(x,y)dy'dy_n\nonumber\\
&\geq&C\int_R^\infty f(u(y_n))y_n^{\alpha/2}\int_{\mathbb{R}^{n-1}\backslash B_R(0)}\frac{1}{|x-y|^n}dy'dy_n.
\label{428}
\end{eqnarray}
When $y_n\rightarrow+\infty$, the behavior of $u(y_n)$ may have two cases:
\begin{description}
\item[~~~(ii)$_1$]$u(y_n)\rightarrow+\infty$;
\item[~~~(ii)$_{2}$]there is a positive constant $c$ such that $u(y_n)\leq c$.
\end{description}
However in both cases, we will be able to derive contradictions in the sequel.

\textbf{(ii)$_1$.} It is easy to see that $u(y_n)\geq C$ for $y_n>R$ and $f(u(y_n))\geq C$ (since $f(t)>0$ is strictly increasing for $t>0$). Set $x=(x',x_n),~y=(y',y_n)\in \mathbb{R}^n_+$, $r^2=|x'-y'|^2$ and $a^2=|x_n-y_n|^2$, then we have form (\ref{428}),
\begin{eqnarray}
+\infty>u(x)&\geq&C\int_R^\infty f(u(y_n))y_n^{\alpha/2}\int_{\mathbb{R}^{n-1}\backslash B_R(0)}\frac{1}{|x-y|^n}dy'dy_n\nonumber\\
&\geq&C\int_R^\infty y_n^{\alpha/2}\int_R^\infty\frac{r^{n-2}}{(r^2+a^2)^\frac{n}{2}}drdy_n\nonumber\\
&\geq&C\int_R^\infty y_n^{\alpha/2}\frac{1}{|x_n-y_n|}\int_{R/a}^\infty\frac{\tau^{n-2}}{(\tau^2+1)^\frac{n}{2}}d\tau dy_n\label{427}\\
&\geq&C\int_R^\infty y_n^{\alpha/2-1}dy_n=\infty,
\label{429}
\end{eqnarray}
which is a contradiction. We derive (\ref{427}) by letting $\tau=\frac{r}{a}$.

\textbf{(ii)$_{2}$.} Since $f(t)/t^{\frac{n+\alpha}{n-\alpha}}$ is non-increasing about $t$ in $(0,+\infty)$ and $u(y_n)\leq c$ for $y_n>R$, it follows that for $y_n>R$,
$$
\frac{f(u(y_n))}{u^{\frac{n+\alpha}{n-\alpha}}(y_n)}\geq\frac{f(c)}{c^{\frac{n+\alpha}{n-\alpha}}}:=C_0,
$$
and from (\ref{428}),
\begin{eqnarray}
+\infty>u(x)&\geq&C\int_R^\infty f(u(y_n))y_n^{\alpha/2}\int_{\mathbb{R}^{n-1}\backslash B_R(0)}\frac{1}{|x-y|^n}dy'dy_n\nonumber\\
&\geq&C\int_R^\infty \frac{f(u(y_n))}{u^{\frac{n+\alpha}{n-\alpha}}(y_n)}u^{\frac{n+\alpha}{n-\alpha}}(y_n)y_n^{\alpha/2}\int_{\mathbb{R}^{n-1}\backslash B_R(0)}\frac{1}{|x-y|^n}dy'dy_n\nonumber\\
&\geq&C\int_R^\infty u^{\frac{n+\alpha}{n-\alpha}}(y_n)y_n^{\alpha/2}\int_{\mathbb{R}^{n-1}\backslash B_R(0)}\frac{1}{|x-y|^n}dy'dy_n\nonumber\\
&\geq&C\int_R^\infty u^{\frac{n+\alpha}{n-\alpha}}(y_n) y_n^{\alpha/2}\int_R^\infty\frac{r^{n-2}}{(r^2+a^2)^\frac{n}{2}}drdy_n\nonumber\\
&\geq&C\int_R^\infty u^{\frac{n+\alpha}{n-\alpha}}(y_n) y_n^{\alpha/2}\frac{1}{|x_n-y_n|}\int_{R/a}^\infty\frac{\tau^{n-2}}{(\tau^2+1)^\frac{n}{2}}d\tau dy_n\nonumber\\
&\geq&C\int_R^\infty u^{\frac{n+\alpha}{n-\alpha}}(y_n)y_n^{\alpha/2-1}dy_n,
\label{42929}
\end{eqnarray}
which implies that there exists a sequence $\{y_n^i\}\rightarrow\infty$ as $i\rightarrow\infty$, such that
\begin{equation}
u^{\frac{n+\alpha}{n-\alpha}}(y_n^i)(y_n^i)^{\alpha/2}\rightarrow0.
\label{42930}
\end{equation}
The rest proof is completely similar to the last part of section 3.2 in \cite{CFY}, here we only need to substitute $u^{\frac{n+\alpha}{n-\alpha}}$ for $u^p$, and we can finally derive a contradiction with (\ref{42930}).

At present, we have proved that $u=u(x_n)$ contradicts the finiteness of $u$. Hence (\ref{03}) possesses no positive solution.

\section{Proofs of Two Claims}
Let us prove Claims 3.1 and 3.2 which are used in Section 3. Without loss of generality, we let $x^0=0$.

\textbf{Proof of Claim 3.1}. Since $v(x)\rightarrow0$ as $|x|\rightarrow\infty$, and $B_\epsilon(0^\lambda)\subset\Sigma_\lambda$ when $\lambda$ is sufficiently negative, we only need to prove $v_\lambda(x)\geq2c_\lambda$ in $B_\epsilon(0^\lambda)\setminus\{0^\lambda\}$, i.e.,
$$v(x)\geq2c_\lambda ~\text{in}~ B_\epsilon(0)\setminus\{0\}.$$
Using
$$
u(x)=\frac{1}{|x|^{n-\alpha}}v(\frac{x}{|x|^2}),
$$
it suffices to prove
$$
u(x)\geq2c_\lambda\frac{1}{|x|^{n-\alpha}},~~~\text{when} ~|x|~\text{sufficiently large.}
$$

Denote
$$
\varphi(x)=c_{n,-\alpha}\int_{\mathbb{R}^n}\frac{\eta(y)f(u(y))}{|x-y|^{n-\alpha}}dy,
$$
where $\eta(y)\in C^\infty(\mathbb{R}^n)$ is a cutoff function and
\begin{equation}
\eta(y)=\left\{\begin{array}{ll}
1, & \qquad y\in B_1(0), \\
0, & \qquad y\notin B_2(0).
\end{array}\right.
\label{jbb2}\end{equation}
Then
\begin{eqnarray}
(-\lap)^{\alpha/2}\varphi(x)&=&c_{n,-\alpha}\int_{\mathbb{R}^n}(-\lap)^{\frac{\alpha}{2}}(\frac{1}{|x-y|^{n-\alpha}})\eta(y)f(u(y))dy\nonumber\\
&=&\eta(x)f(u(x))
\end{eqnarray}
and
\begin{eqnarray}
(-\lap)^{\alpha/2}(u-\varphi)&=&f(u)-\eta f(u)\nonumber\\
&=&f(u)(1-\eta)\nonumber\\
&\geq&0,~~~\text{in}~\mathbb{R}^n.
\label{APp1}
\end{eqnarray}
For any $x\in B_R^c(0)$, we have
$$
c_{n,-\alpha}\int_{B_1(0)}\frac{f(u(y))}{|x-y|^{n-\alpha}}dy\leq\varphi(x)\leq c_{n,-\alpha}\int_{B_2(0)}\frac{f(u(y))}{|x-y|^{n-\alpha}}dy.
$$
Because $u\in C_{loc}^{1,1}(\mathbb{R}^n)$, $f>0$ is strictly increasing for $t>0$, it yields that for any $x\in B_R^c(0)$ and $R$ large enough, there exist two positive constants $C_1$ and $C_2$ such that
\begin{equation}
\frac{C_1}{|x|^{n-\alpha}}\leq\varphi(x)\leq\frac{C_2}{|x|^{n-\alpha}}\leq\frac{C_2}{R^{n-\alpha}}.
\label{App3}
\end{equation}
Combining (\ref{APp1}) and (\ref{App3}), we get
\begin{equation}
\left\{\begin{array}{ll}
(-\lap)^{\alpha/2}(u-\varphi+\frac{C_2}{R^{n-\alpha}})\geq0, & \qquad \text{in}~ B_R(0), \\
u-\varphi+\frac{C_2}{R^{n-\alpha}}\geq0, & \qquad \text{in}~ B_R^c(0)
\end{array}\right.
\label{jbb2}\end{equation}
and
$$
u-\varphi+\frac{C_2}{R^{n-\alpha}}\geq0,~~\text{in}~ B_R(0),
$$
by the \textsl{maximum principle} in \cite{Si}. Letting $R\rightarrow\infty$, it derives
$$
u\geq\varphi,~~~x\in\mathbb{R}^n,
$$
and by (\ref{App3}),
$$
u(x)\geq\frac{2c_\lambda}{|x|^{n-\alpha}},~~~\text{for}~|x|~\text{sufficiently large}.
$$
This complete the proof of Claim 3.1.

Before proving Claim 3.2, we first narrate
\begin{proposition}
Assume that $u\in\mathcal{L}_\alpha\cap C_{loc}^{1,1}(\mathbb{R}^n)$ is a locally bounded positive solution to
\begin{equation}
(-\lap)^{\alpha/2} u(x)=f(u), \qquad x\in\mathbb{R}^n,
\label{A423}
\end{equation}
then it is also a solution to
\begin{equation}
u(x)=\int_{\mathbb{R}^n}\frac{C}{|x-y|^{n-\alpha}}f(u(y))dy.
\label{A424}
\end{equation}
And vice versa.
\label{pro5.1}
\end{proposition}

The proof of Proposition \ref{pro5.1} is entirely similar to the proof of Theorem 2.1 in \cite{ZC} and we only need to guarantee that: ($i$) for any positive constant $c$, $f(c)>0$; ($ii$) $f(u(x))$ is locally bounded on $\mathbb{R}^n$. These can be easily acquired by the condition on $f(t)$. Here we omit the details.

\textbf{Proof of Claim 3.2.} If
$$
w_{\lambda_0}(x)\not\equiv0,~~x\in\Sigma_{\lambda_0}\setminus\{0^{\lambda_0}\},
$$
then there exists a point $x^1\in\Sigma_{\lambda_0}\setminus\{0^{\lambda_0}\}$ such that
$$
w_{\lambda_0}(x^1)>0,
$$
and so there is a small positive $\delta$ such that
\begin{equation}
w_{\lambda_0}(x)>0,~~x\in B_{\delta}(x^1),
\end{equation}
and
\begin{equation}
w_{\lambda_0}(x)\geq c_1,~~x\in B_{\delta/2}(x^1).
\label{AA1}
\end{equation}
From Proposition \ref{pro5.1} and (\ref{3002}) (please note that we already let $x^0=0$ for simplicity at the beginning of this section), we have
\begin{equation}
(-\lap)^{\alpha/2} v(x)=\int_{\mathbb{R}^n}\frac{C}{|x-y|^{n-\alpha}}\frac{1}{|y|^{n+\alpha}}f(|y|^{n-\alpha}v(y))dy,\quad x\in\mathbb{R}^n\setminus\{0\}.
\label{300200}
\end{equation}
Through an elementary computation, we have that for any $x\in\mathbb{R}^n\setminus\{0^{\lambda_0}\}$,
\begin{equation}
w_{\lambda_0}(x)=C\int_{\Sigma_{\lambda_0}}(\frac{1}{|x-y|^{n-\alpha}}-\frac{1}{|x-y^{\lambda_0}|^{n-\alpha}})(\frac{f(|y^{\lambda_0}|^{n-\alpha}v(y^{\lambda_0}))}{|y^{\lambda_0}|^{n+\alpha}}-\frac{f(|y|^{n-\alpha}v(y))}{|y|^{n+\alpha}})dy.
\end{equation}
It follows from the monotonicities of $f(t)$ and $f(t)/t^{\frac{n+\alpha}{n-\alpha}}$ that for any $x\in B_\epsilon(0^{\lambda_0})\setminus\{0^{\lambda_0}\}$,
\begin{eqnarray}
w_{\lambda_0}(x)&=&C\int_{\Sigma_{\lambda_0}}(\frac{1}{|x-y|^{n-\alpha}}-\frac{1}{|x-y^{\lambda_0}|^{n-\alpha}})(\frac{f(|y^{\lambda_0}|^{n-\alpha}v(y^{\lambda_0}))}{|y^{\lambda_0}|^{n+\alpha}}-\frac{f(|y|^{n-\alpha}v(y))}{|y|^{n+\alpha}})dy\nonumber\\
&\geq&C\int_{\Sigma_{\lambda_0}}(\frac{1}{|x-y|^{n-\alpha}}-\frac{1}{|x-y^{\lambda_0}|^{n-\alpha}})(\frac{f(|y|^{n-\alpha}v(y^{\lambda_0}))-f(|y|^{n-\alpha}v(y))}{|y|^{n+\alpha}})dy\nonumber\\
&\geq&C\int_{B_{\delta/2}(x^1)}(\frac{1}{|x-y|^{n-\alpha}}-\frac{1}{|x-y^{\lambda_0}|^{n-\alpha}})(\frac{f(|y|^{n-\alpha}v(y^{\lambda_0}))-f(|y|^{n-\alpha}v(y))}{|y|^{n+\alpha}})dy.\nonumber\\
\label{App7}
\end{eqnarray}
For any $y\in\overline{B_{\delta/2}(x^1)}$,
$$|y|^{n-\alpha}v(y^{\lambda_0})-|y|^{n-\alpha}v(y)=|y|^{n-\alpha}w_{\lambda_0}(y),$$
one knows by the continuity of $v$ and (\ref{AA1}) that
$$c_2\leq|y|^{n-\alpha}v(y^{\lambda_0})-|y|^{n-\alpha}v(y)\leq c_3$$
for some constants $c_2, c_3>0$, $c_2<c_3$. Since $f$ is strictly increasing, there exists $c_4,~c_5>0$, $c_4<c_5$ such that for any $x\in\overline{B_{\delta/2}(x^1)}$,
$$c_4\leq f(|y|^{n-\alpha}v(y^{\lambda_0}))-f(|y|^{n-\alpha}v(y))\leq c_5.$$
Hence
\begin{eqnarray*}
f(|y|^{n-\alpha}v(y^{\lambda_0}))-f(|y|^{n-\alpha}v(y))&\geq& C_1(|y|^{n-\alpha}v(y^{\lambda_0})-|y|^{n-\alpha}v(y))\\
&=&C_1(|y|^{n-\alpha}w_{\lambda_0}(y)),
\end{eqnarray*}
for some $C_1>0$. It shows by (\ref{AA1}) and (\ref{App7}) that
\begin{eqnarray}
w_{\lambda_0}(x)&\geq&C\int_{B_{\delta/2}(x^1)}(\frac{1}{|x-y|^{n-\alpha}}-\frac{1}{|x-y^{\lambda_0}|^{n-\alpha}})(\frac{f(|y|^{n-\alpha}v(y^{\lambda_0}))-f(|y|^{n-\alpha}v(y))}{|y|^{n+\alpha}})dy\nonumber\\
&\geq&C\int_{B_{\delta/2}(x^1)}(\frac{1}{|x-y|^{n-\alpha}}-\frac{1}{|x-y^{\lambda_0}|^{n-\alpha}})\frac{C_1}{|y|^{2\alpha}}w_{\lambda_0}(y)dy\nonumber\\
&\geq&C\int_{B_{\delta/2}(x^1)}C_1dy\nonumber\\
&:=&2c_0
\label{App8}
\end{eqnarray}
and Claim 3.2 is proved.




\end{document}